\documentclass{article}

\usepackage{proof}
\usepackage{latexsym}
\usepackage{amssymb}

\usepackage[all]{xy}

\newcommand{\blem}{\begin{lemma}}
\newcommand{\elem}{\end{lemma}}
\newcommand{\bth}{\begin{theorem}}
\newcommand{\ethm}{\end{theorem}}
\newcommand{\benu}{\begin{enumerate}}
\newcommand{\eenu}{\end{enumerate}}
\newcommand{\bdes}{\begin{description}}
\newcommand{\edes}{\end{description}}
\newcommand{\bdf}{\begin{definition}}
\newcommand{\edf}{\end{definition}}
\newcommand{\bcor}{\begin{cor}}
\newcommand{\ecor}{\end{cor}}
\newcommand{\bprp}{\begin{proposition}}
\newcommand{\eprp}{\end{proposition}}
\newcommand{\bmlem}{\begin{mlemma}}
\newcommand{\emlem}{\end{mlemma}}
\newcommand{\bclm}{\begin{claim}}
\newcommand{\eclm}{\end{claim}}
\newcommand{\bprf}{{\bf Proof}.\hspace{2mm}}

\newcommand{\eprf}{\hspace*{\fill} $\Box$}

\newcommand{\beqn}{\begin{equation}}
\newcommand{\eeqn}{\end{equation}}
\newcommand{\beqnarr}{\begin{eqnarray}}
\newcommand{\eeqnarr}{\end{eqnarray}}
\newcommand{\beqnarrs}{\begin{eqnarray*}}
\newcommand{\eeqnarrs}{\end{eqnarray*}}

\newcommand{\spand}{\,\&\,}

\newcommand{\Natural}{\mathbb{N}}

\newtheorem{theorem}{Theorem}[section]
\newtheorem{definition}[theorem]{Definition}
\newtheorem{proposition}[theorem]{Proposition}
\newtheorem{lemma}[theorem]{Lemma}
\newtheorem{cor}[theorem]{Corollary}
\newtheorem{mlemma}[theorem]{Main Lemma}
\newtheorem{claim}[theorem]{Claim}

\newcommand{\alp}{\alpha}

\newcommand{\veps}{\varepsilon}

\newcommand{\ome}{\omega}

\newcommand{\bet}{\beta}
\newcommand{\gam}{\gamma}

\newcommand{\Sig}{\Sigma}

\newcommand{\fal}{\forall}
\newcommand{\exi}{\exists}

\newcommand{\Rarw }{\Rightarrow}

\newcommand{\Lrarw}{\Leftrightarrow}

\newcommand{\brem}{\begin{remark}}
\newcommand{\erem}{\end{remark}}
\newtheorem{remark}[theorem]{Remark}

\title{Grzegorczyk sequence
}
\author{
Toshiyasu Arai
\\
Graduate School of Science,
Chiba University
\\
1-33, Yayoi-cho, Inage-ku,
Chiba, 263-8522, JAPAN
\\
tosarai@faculty.chiba-u.jp
}
\date{}

\begin{document}
\maketitle
\begin{abstract}
Natural numbers are represented by Grzegorczyk functions.
The representation is implicit in the technique of H. Friedman\cite{Simpson85}.
An iterated base-shift in the representation with subtracting $1$ yields a sequence, \textit{Grzegorczyk sequence}.
It is shown that the termination of the sequence is independent from the first order arithmetic {\sf PA}.
We follow M. Rathjen\cite{Rathjen15} in the proof of the independence. 
\end{abstract}

\section{Grzegorczyk representation of natural numbers and Grzegorczyk sequences}
Let us represent natural numbers in terms of Grzegorczyk functions $F_{n}$.
The representation is implicit in the technique of H. Friedman\cite{Simpson85}.

\bdf\label{df:Grzegorczyk}
{\rm
\benu
\item
$f^{(x)}$ denotes the $x$-th iterate of a unary function $f$ on $\Natural$, defined by
$f^{(0)}(y)=y$ and $f^{(x+1)}(y)=f(f^{(x)}(y))$.
\item
The \textit{$n$-th Grzegorczyk function}
$F_{n}$ denote the function on integers defined recursively on $n$ as follows.
$F_{0}(x)=x+1$ and $F_{n+1}(x)=F_{n}^{(x)}(x)$.
\eenu
}
\edf

\bprp\label{prp:Grzegorczyk}
\benu
\item\label{prp:Grzegorczyk.1}
$F_{n}(x)>x$ for $x>0$.
\item\label{prp:Grzegorczyk.2}
$F_{n}^{(y)}(x)\geq x$.
\item\label{prp:Grzegorczyk.3}
$(x,y)\mapsto F_{n}^{(y)}(x)$ is strictly increasing both in $x>0$ and $y$.
\item\label{prp:Grzegorczyk.4}
$F_{n+1}(x)> F_{n}(x)$ for $x\geq 2$.
\eenu
\eprp

\blem\label{lem:Grzegorczyk}
For each primitive recursive function $f(\vec{x})$, there exists an $n$ for which
$f(\vec{x})\leq F_{n}(\max\{\vec{x},2\})$ holds for any $\vec{x}=x_{1},\ldots,x_{k}$.
\elem

\bdf\label{df:Frep}
{\rm
Let $k_{1}\geq 2$.
$F_{n}$ denotes the $n$-th Grzegorczyk function.
We define \textit{$F$-representation} 
of natural numbers $x$ \textit{with base $k_{1}$}
by induction on $x\dot{-}k_{1}$
as follows.

First if $x<k_{1}$, then $x$ itself is the $F$-representation of $x$ with base $k_{1}$.
Assume $k_{1}\leq x$.
Let $x_{1}$ be the least natural number such that $k_{1}\leq x<F_{x_{1}+1}(k_{1})$, and
$i_{1}$ the maximal number such that 
$F_{x_{1}}^{(i_{1})}(k_{1})\leq x<F_{x_{1}}^{(k_{1})}(k_{1})=F_{x_{1}+1}(k_{1})$.
If $x_{1}=0$, then $[(0,i_{1})]_{k_{1}}$ with $0=x_{1}$ is the $F$-representation of $x$.
Otherwise let $x_{2}<x_{1}$ be the least natural number such that
$k_{2}\leq x<F_{x_{2}+1}(k_{2})$ for $k_{2}=F_{x_{1}}^{(i_{1})}(k_{1})\leq x<F_{x_{1}}(k_{2})$.
If $x=k_{2}$, then $[(x_{1},i_{1})]_{k_{1}}$ is the $F$-representation of $x$ with base $k_{1}$.
Otherwise
$[(x_{1},i_{1}),(x_{2},i_{2}),\ldots,(x_{\ell},i_{\ell})]_{k_{1}}$ is the 
$F$-representation of $x$ with base $k_{1}$, 
where $[(x_{2},i_{2}),\ldots,(x_{\ell},i_{\ell})]_{k_{2}}\neq[(0,0)]_{k_{2}}$ is the $F$-representation of $x$ with base $k_{2}$.
}
\edf

We define \textit{Grzegorczyk sequences} $\{z_{k}\}_{k}$ of natural numbers $z$ as follows.
$z_{k+1}=z_{k}\dot{-}1$ if $z_{k}<2+k$.
Let $z_{k}\geq 2+k$.
Let $x=[(x_{1},i_{1}),(x_{2},i_{2}),\ldots,(x_{\ell},i_{\ell})]_{k}$ be the $F$-representation of 
$x\geq k$ with base $k$.
Then the shift of base $k$ to $m$ in the representation as well as the representations of $x_{i}$
hereditarily in the left is defined recursively by
$x[k:=m]=
[(x_{1}[k:=m],i_{1}),(x_{2}[k:=m],i_{2}),\ldots,(x_{\ell}[k:=m],i_{\ell})]_{m}$.
Then
$z_{k+1}=z_{k}[2+k:=3+k]-1$.

Since the function $(n,x)\mapsto F_{n}(x)$ is provably computable in {\sf PA}
(but not primitive recursive), so is the function $(x,k,m)\mapsto x[k:=m]$.

\bth\label{th:GrzegorczykseqPA}
${\sf PA}$ does not prove the true $\Pi_{2}$-statement $\fal z \exi k[z_{k}=0]$
for the Grzegorczyk sequence $\{z_{k}\}_{k}$ of $z$.
\end{theorem}

\section{Slowly well-foundedness}
\bdf\label{df:normepsilon}
{\rm
We define integers $C(\alp)$ for ordinal terms $\alp<\veps_{0}$.
$C(0)=0$.
In the following let $\alp=\ome^{\alp_{1}}n_{1}+\cdots+\ome^{\alp_{k}}n_{k}$
with $\alp_{1}>\cdots>\alp_{k}$ and $0<n_{1},\ldots,n_{k}<\ome$ for $k>0$.
$C(\alp)=\max\{C(\alp_{\ell}), n_{\ell}:1\leq\ell\leq k\}$.
$C(\alp)$ is the maximal coefficient $n_{\ell}$ in the ordinal term $\alp$.
}
\edf

In the $F$-representation of Definition \ref{df:Frep},
for the case $k_{1}\leq x$,
note that $k_{1}<F_{n}(k_{1})<F_{n+1}(k_{1})$ holds by $k_{1}\geq 2$ and Proposition \ref{prp:Grzegorczyk}.\ref{prp:Grzegorczyk.1}.
When $x_{1}>0$, $i_{1}>0$ holds by the minimality of $x_{1}$.
Hence $k_{1}<F_{x_{1}}^{(i_{1})}(k_{1})=k_{2}$.

It is clear that for $x\geq k_{1}$, $[(0,0)]_{k_{1}}$ is the $F$-representation of $x$ iff $x=k_{1}$.
Let $[(x_{1},i_{1}),(x_{2},i_{2}),\ldots,(x_{\ell},i_{\ell})]_{k_{1}}\neq[(0,0)]_{k_{1}}$
 be the $F$-representation of
$x$ with base $k_{1}$ for $2\leq k_{1}< x$.
Proposition \ref{prp:Grzegorczyk}.\ref{prp:Grzegorczyk.4} yields $x_{1}<F_{x_{1}}(k_{1})\leq x$.
Hence $x>x_{1}>x_{2}>\cdots>x_{\ell}\geq 0$,
$\fal p\leq\ell(i_{p}>0)$.

Let $k_{\ell}\, (1\leq p\leq\ell+1)$ denote the integers recursively defined as follows.
$k_{1}$ is the given number, and
$k_{p+1}=F_{x_{p}}^{(i_{p})}(k_{p})$.
Then
$i_{p}<k_{p}$, and $k_{1}<k_{2}<\cdots<k_{\ell}<k_{\ell+1}=x$, i.e.,
$x=F_{x_{\ell}}^{(i_{\ell})}(\cdots(F_{x_{1}}^{(i_{1})}(k_{1}))\cdots)$.
Moreover $\fal p\leq \ell(i_{p}<k_{p}<x)$.
Hence we obtain in the $F$-representation $[(x_{1},i_{1}),(x_{2},i_{2}),\ldots,(x_{\ell},i_{\ell})]_{k_{1}}$
of $x\geq k_{1}\geq 2$
\beqn\label{eq:Frepless}
x>x_{1}>x_{2}>\cdots>x_{\ell}\geq 0 \spand \fal p\leq \ell(i_{p}<x)
\eeqn
Note that the function assigning $x\geq k$ to its $F$-representation
\\
$[(x_{1},i_{1}),(x_{2},i_{2}),\ldots,(x_{\ell},i_{\ell})]_{k_{1}}$
is elementary recursive since the ternary relation $R=\{(n,x,y)\in\Natural^{3}:F_{n}(x)=y\}$
is elementary recursive.
This is seen as follows.
By Proposition \ref{prp:Grzegorczyk}.\ref{prp:Grzegorczyk.1} we see that
$(n,x,y)\in R$ iff there exists a matrix $A=(a_{ij})$ such that its size is at most $n\times y$,
and its entry $a_{ij}\leq\max\{x,y\}$ is a number $F_{m}(y_{ij})$ appearing in the computation of
$F_{n-1}^{(p)}(x)\,(p\leq x)$.

Thus $x$ is represented by
\[
[(x_{1},i_{1}),(x_{2},i_{2}),\ldots,(x_{\ell},i_{\ell})]_{k_{1}}:=
F_{x_{\ell}}^{(i_{\ell})}(\cdots(F_{x_{1}}^{(i_{1})}(k_{1}))\cdots)=x
.\]

Let
$[(x_{1},i_{1}),(x_{2},i_{2}),\ldots,(x_{\ell},i_{\ell})]_{k_{1}}$ and
$[(y_{1},j_{1}),(y_{2},j_{2}),\ldots,(y_{m},j_{m})]_{k_{1}}$ be the $F$-representations of 
$x\geq k_{1}$ and $y\geq k_{1}$, resp.
Then we see the following.
\beqn\label{eq:Frep}
x<y \Lrarw 
\left((x_{1},i_{1}),(x_{2},i_{2}),\ldots,(x_{\ell},i_{\ell})\right)
\prec
\left((y_{1},j_{1}),(y_{2},j_{2}),\ldots,(y_{m},j_{m})\right)
\eeqn
where for the lexicographic order $<_{lx}$ on pairs,
\[
\left((x_{1},i_{1}),(x_{2},i_{2}),\ldots,(x_{\ell},i_{\ell})\right)
\prec
\left((y_{1},j_{1}),(y_{2},j_{2}),\ldots,(y_{m},j_{m})\right)
\]
iff either
$\exi p\leq\min\{\ell,m\}[\fal q<p\{(x_{q},i_{q})=(y_{q},j_{q})\}\land (x_{p},i_{p})<_{lx}(y_{p},j_{p})]$
or $\ell<m\land \fal q\leq\ell\{(x_{q},i_{q})=(y_{q},j_{q})\}$.

Let $k_{1}\leq x<F_{n}(k_{1})$, and
$[(x_{1},i_{1}),(x_{2},i_{2}),\ldots,(x_{\ell},i_{\ell})]_{k_{1}}$ the $F$-representation of $x$
with base $k_{1}$, where $n> x_{1}>x_{2}>\cdots>x_{\ell}\geq 0$.
$s(x;n,k_{1}):=(j_{1},j_{2},\ldots,j_{n})$ denotes the sequence of natural numbers $j_{q}$,
which is obtained from the sequence $(i_{1},i_{2},\ldots,i_{\ell})$ by filling the gaps with zeros,
i.e., $j_{x_{p}}=i_{p}$ and $j_{q}=0$ else.
Moreover let $(m_{1},m_{2},\ldots,m_{n})$ be the sequence defined recursively by
$m_{1}=k_{1}$ and $m_{q+1}=F_{n-q}^{(j_{q})}(m_{q})$.

Let $y_{1},y_{2}$ be integers such that $k_{1}\leq y_{1},y_{2}<F_{n}(k_{1})$ for $n>0$.
For simplicity let us write $s(x)$ for $s(x;n,k_{1})$.
We see from Proposition \ref{prp:Grzegorczyk}.\ref{prp:Grzegorczyk.1} that for $k_{1}\geq 2$
\beqn\label{eq:slows}
y_{1}<y_{2} \Lrarw s(y_{1})<_{lx}s(y_{2})
\eeqn
where $<_{lx}$ denotes the lexicographic order on $n$-tuples.

Now let us flip over the integers $i_{\ell}$.
Namely let ${}^{t}s(x)=(m_{1}-j_{1},m_{2}-j_{2},\ldots,m_{n}-j_{n})$ for
$s(x)=(j_{1},j_{2},\ldots,j_{n})$.
${}^{t}s(x)$ is an $n$-tuple of positive integers $m_{q}-j_{q}\leq x$ such that
\[
y_{1}<y_{2}\Lrarw {}^{t}s(y_{2})<_{lx} {}^{t}s(y_{1})
\]
We obtain a descending chain $\{g_{n}(k_{1},x)\}_{x<F_{n}(k_{1})}$ of ordinals
$g_{n}(k_{1},x)<\ome^{n+1}$ for $n>0$ as follows.
When $k_{1}\leq x<F_{n}(k_{1})$, let 
\beqn\label{eq:gn}
g_{n}(k_{1},x)=\ome^{n-1}(m_{1}-j_{1})+\ome^{n-2}(m_{2}-j_{2})+\cdots+\ome^{0}(m_{n}-j_{n})
\eeqn
for $s(x;n,k_{1})=(j_{1},j_{2},\ldots,j_{n})$.
When $x<k_{1}$, let
$g_{n}(k_{1},x)=\ome^{n}(k_{1}-x)$.

Let us extend $g_{n}$ by letting $g_{n}(k_{1},x)=0$ for $x\geq F_{n}(k_{1})$, and
$g_{0}(k_{1},x)=(k_{1}+1)\dot{-}x$.
Then it is easy to see that for each $n$ there exists a constant $c_{n}$ for which 
$C(g_{n}(k_{1},x))\leq\max\{n,k_{1}+1,x\}$ holds
for any $k_{1},x$.

\blem\label{lem:slowdown}{\rm (H. Friedman\cite{Simpson85})}\\
For each unary primitive recursive function $f$, there exist $n, c$ and a primitive recursive function
$g:\Natural^{2}\to\ome^{n}$ such that
$x<f(k) \Rarw g(k,x+1)<g(k,x)$, and
$C(g(k,x))\leq\max\{c,k+1,x\}$ for any $k,x$.
\elem
\bprf
Pick an $n$ such that $f(k)\leq F_{n}(\max\{2,k\})$ by Lemma \ref{lem:Grzegorczyk}.
\eprf

\bth\label{th:slowdown}
Let $\veps_{0}>\alp_{0}>\alp_{1}>\cdots$ be a primitive recursive descending chain of ordinal terms $\alp_{k}$.
We can find a primitive recursive descending chain $\{\gam_{i}\}_{i}$ of ordinal terms such that
$C(\gam_{i})\leq i+1$.
\end{theorem}
\bprf
Consider the primitive recursive functions 
$f_{C}(k)=C(\alp_{k+1})$.
By Lemma \ref{lem:slowdown}, 
let $h:\Natural^{2}\to\ome^{n}$ be a primitive recursive function and $c$ a constant such that
$x<C(\alp_{k+1}) \Rarw h(k,x+1)<h(k,x)$, $C(h(k,x))\leq\max\{c,k+1,x\}$.

Llet $\ell=\max\{c,C(\alp_{0})\}$.
For $i\geq \ell$, let $k,x$ be numbers such that
$i=C(\alp_{0})+\cdots+C(\alp_{k})+x$ and $x<C(\alp_{k+1})$.
Then set 
\[
\gam_{i}  =  \ome^{\ome}\alp_{k}+h(k,x)
\]
Note here that $C(\alp_{k})>0$ for $\alp_{k}>0$.

Let $N$ be a number such that $\ome_{N-\ell}>\ome^{\ome}\alp_{0}$.
For $i<\ell$ let 
$\gam_{i}=\ome_{N-i}$.
Then the primitive recursive sequence $\{\gam_{i}\}_{i}$ is descending, and $C(\gam_{i})\leq i+1$.
\eprf

\bcor\label{cor:slowdown}
Over {\sf PRA}, the 1-consistency $\mbox{{\rm RFN}}_{\Sig_{1}}({\sf PA})$ of {\sf PA}
is equivalent to the fact that
there is no primitive recursive and 
infinite descending chain $\{\gam_{i}\}_{i}$ of ordinals $\gam_{i}<\veps_{0}$
such that $C(\gam_{i})\leq i+1$.
\ecor
\bprf
This is seen from the fact that $\mbox{{\rm RFN}}_{\Sig_{1}}({\sf PA})$ is equivalent to the fact that 
there is no primitive recursive and 
infinite descending chain of ordinals$<\veps_{0}$.
\eprf

\section{Proof of Theorem \ref{th:GrzegorczykseqPA}}

For $x\geq k\geq 2$,
let $x=[(x_{1},i_{1}),(x_{2},i_{2}),\ldots,(x_{\ell},i_{\ell})]_{k}$
be the $F$-representation of $x$ with base $k$.
Then an ordinal $o_{k}(x)<\veps_{0}$ is associated recursively by
\[
o_{k}(x)=\ome^{o_{k}(x_{1})}i_{1}+\ome^{o_{k}(x_{2})}i_{2}+\cdots+\ome^{o_{k}(x_{\ell})}i_{\ell}
\]
From (\ref{eq:Frep}) we see that for $x,y\geq k\geq 2$,
$o_{k}(x)$ is in Cantor normal form with base $\ome$ and
\beqn\label{eq:Disoleft}
x<y\Lrarw o_{k}(x)<o_{k}(y)
\eeqn
Let 
\[
D_{k}:=\{o_{k}(x):x\geq k\}
\]
and $L_{k}(\alp)$ be the inverse of $o_{k}$ for $\alp\in D_{k}$, $L_{k}(o_{k}(x))=x$.

Let us prove Theorem \ref{th:GrzegorczykseqPA}.
First we show the termination of Grzegorczyk sequences $\{z_{k}\}_{k}$ of $z\geq 2$.
Suppose $\fal k(z_{k}>0)$. We have $\fal k(z_{k}\geq 2+k)$, for
if $z_{k}=m<2+k$, then $z_{k+m}=0$.
Hence we obtain $\fal k(0<o_{2+k}(z_{k})\in D_{2+k})$.

From (\ref{eq:Disoleft}) we see that $\ome$ is the order type of the set $D_{k}$ of ordinals$<\veps_{0}$ when $k\geq 2$.
For $0<\alp\in D_{k}$, let $Q_{k}\alp=\max\{\bet<\alp: \bet\in D_{k}\}$.
Then it is easy to see that
$z_{k+1}=L_{3+k}Q_{3+k}o_{2+k}(z_{k})$ from $o_{2+k}(z_{k})\in D_{2+k}\subset D_{3+k}$.
Therefore $\{o_{2+k}(z_{k})\}_{k}$ would be an infinite descending chain of ordinals$<\veps_{0}$.
Thus we show the fact that each Grzegorczyk sequences $\{z_{k}\}_{k}$ of $z$ eventually terminates.

Next let us show the independence following M. Rathjen\cite{Rathjen15}.
For the independence of the fact from ${\sf PA}$, it suffices to show that
the fact implies the $1$-consistency of ${\sf PA}$ over {\sf PA}.
Argue in {\sf PA}.
Suppose that ${\sf PA}$ is $1$-inconsistent.
By Corollary \ref{cor:slowdown}
let $\{\gam_{k}\}_{k}$ be a primitive recursive and infinite descending chain of
ordinals $\gam_{i}<\veps_{0}$ such that $C(\gam_{k})\leq 1+k$.
Then $\gam_{k}\in D_{2+k}$ since $k+2<F_{x}^{(i)}(k+2)$ for $i>0$.
Hence $\{\gam_{k}\}_{k}$ would be a primitive recursive and infinite descending chain of
ordinals$<\veps_{0}$ such that $\gam_{k}\in D_{2+k}$.
Let $v_{k}=L_{2+k}\gam_{k}$, and $\{z_{k}\}_{k}$ be the Grzegorczyk sequence of $z_{0}=v_{0}$.
We show that $\fal k(v_{k}\leq z_{k})$ by induction on $k$.
Suppose $v_{k}\leq z_{k}$. Then $\gam_{k}=o_{2+k}(v_{k})\leq o_{2+k}(z_{k})$.
From $\gam_{k},\gam_{k+1}\in D_{3+k}$ and $\gam_{k+1}<\gam_{k}\leq o_{2+k}(z_{k})$ we obtain 
$\gam_{k+1}\leq Q_{3+k}o_{2+k}(z_{k})$.
Hence $v_{k+1}=L_{3+k}\gam_{k+1}\leq L_{3+k}Q_{3+k}o_{2+k}(z_{k})=z_{k+1}$.
Thus $\fal k(v_{k}\leq z_{k})$ is shown.
Since $\gam_{k}> 0$, we have $v_{k}>0$.
Thus $\{z_{k}\}_{k}$ would be a non-terminating Grzegorczyk sequence of $z_{0}$.
\eprf
\\

\noindent
The \textit{total representation}  $T(x)_{k_{1}}$ of $x$ \textit{with base} $k_{1}\geq 2$
represents $x$  hereditarily in the both components $x_{p},i_{p}$.
These are obtained by modifying Definition \ref{df:Frep} of $F$-representations as follows.

Let $k\geq 2$.
First if $x<k$, then $T(x)_{k}=x$ itself is 
the hereditary $F$-representations of $x$ with base $k$.
Assume $k\leq x$, and let $[(x_{1},i_{1}),(x_{2},i_{2}),\ldots,(x_{\ell},i_{\ell})]_{k}$
be the $F$-representation of $x$ with base $k$.
Define hereditary representations recursively as follows, cf.\,(\ref{eq:Frepless}).
\\
$T(x)_{k}:=[(T(x_{1})_{k},T(i_{1})_{k}),(T(x_{2})_{k},T(i_{2})_{k}),\ldots,(T(x_{\ell})_{k},T(i_{\ell})_{k})]_{k}$,
\\
where $T(x_{p})_{k}, T(i_{p})_{k}$ are the total representations of $x_{p},i_{p}$ with base $k$.

We define \textit{hereditarily Grzegorczyk sequences} $\{w_{k}\}_{k}$ of natural numbers $w$ as follows.
$w_{k+1}=w_{k}\dot{-}1$ if $w_{k}<2+k$.
Let $w_{k}\geq 2+k$.

Let $x=[(x_{1},i_{1}),(x_{2},i_{2}),\ldots,(x_{\ell},i_{\ell})]_{k}$ be the hereditary $F$-representation of 
$x\geq k$ with base $k$.
Then the shift of base $k$ to $m$ in the \textit{both} representations of $x_{p}, i_{p}$
hereditarily is defined recursively by
$x[k:=m]=
[(x_{1}[k:=m],i_{1}[k:=m]),(x_{2}[k:=m],i_{2}[k:=m]),\ldots,(x_{\ell}[k:=m],i_{\ell}[k:=m])]_{m}$.
Then
$w_{k+1}=w_{k}[2+k:=3+k]-1$.
\\

\noindent
{\bf Open problem}.
How strong the true statement $\fal w\exi k[w_{k}=0]$ for the hereditary Grzegorczyk sequence $\{w_{k}\}$ of $w$?
Specifically
does ${\sf ID}_{1}$ prove the $\Pi_{2}$-statement?

\end{document}